\documentclass[12pt,a4paper]{article}

\usepackage[authordate-trad,backend=biber]{biblatex-chicago}
\bibliography{storage_refs_all}

\usepackage{amsmath,amsthm,amsfonts}
\usepackage{graphicx}
\usepackage{color}
\usepackage{setspace}
\usepackage{paralist}
\usepackage{hyperref}
\hypersetup{%
plainpages=false,
colorlinks,urlcolor=blue,linkcolor=blue,citecolor=blue,
pdftitle={The integration of variable generation and storage into
  electricity capacity markets},
pdfauthor={Stan Zachary, Amy Wilson and Chris Dent},
pdfsubject={},
pdfpagemode={UseOutlines},
pdfpagelayout={OneColumn},
pdfstartview={Dubhe}
}

\newtheorem{proposition}{Proposition}
\theoremstyle{remark}

\renewcommand{\Pr}{\mathbf{P}}
\newcommand{\E}{\mathbf{E}}

\newcommand{\R}{R}
\newcommand{\Ri}{R+i}
\newcommand{\Rj}{R+j}
\newcommand{\Rij}{R+i+j}
\renewcommand{\S}{S}

\newcommand{\lole}{\ensuremath{\mathrm{LOLE}}}
\newcommand{\eeu}{\ensuremath{\mathrm{EEU}}}

\newcommand{\cone}{\ensuremath{\mathrm{CONE}}}
\newcommand{\voll}{\ensuremath{\mathrm{VOLL}}}
\newcommand{\efc}[2]{\ensuremath{\mathit{efc}_{#2}(#1)}}
\newcommand{\elcc}[2]{\ensuremath{\mathit{elcc}_{#2}(#1)}}

\parskip\smallskipamount
\parskip \smallskipamount

\title{The integration of variable generation and storage into
  electricity capacity markets}

\author{Stan Zachary\footnote{Corresponding author.  Department of
    Mathematics, University of Edinburgh, Edinburgh, EH9 3FD, UK.
    Email: s.zachary@gmail.com.}, 
  Amy Wilson\footnote{University of Edinburgh.},
  and Chris Dent\footnotemark[2]}
\date{\today}

\begin{document}

\maketitle

\begin{abstract}\noindent
  We show how to value both variable generation and energy storage to
  enable them to be integrated fairly and optimally into electricity
  capacity markets.  We develop theory based on balancing expected
  energy unserved against costs of capacity procurement, and in which
  the optimal procurement is that necessary to meet an appropriate
  reliability standard.  For conventional generation the theory
  reduces to that already in common use.  Further the valuation of
  both variable generation and storage coincides with the traditional
  risk-based approach based on equivalent firm capacity.  The
  determination of the equivalent firm capacity of storage requires
  particular care; this is due both to the flexibility with which
  storage added to an existing system may be scheduled, and also
  because, when \emph{any} resource is added to an existing system,
  storage already within that system may be flexibly rescheduled.  We
  illustrate the theory with an example based on the GB system.

\end{abstract}

\newpage

\section{Introduction}
\label{sec:introduction}

In order to ensure the adequacy of electricity supplies many systems,
including those of Great Britain and other European countries and
several North American regions, now provide a \emph{capacity
  market}---see, e.g.,
\cite{Helm,Newbery2015,MoyeMeyn2018,HolmbergRitz,NGCM18,ISONE}.
(Others, e.g.\ those of Australia and Texas, continue to rely on
energy-only markets.)
So as to operate such a market both fairly and optimally it is
necessary to value appropriately (i.e.\ mathematically correctly) the contributions of the individual
capacity providers, whether they provide conventional generation,
variable generation, or storage.  Present approaches to capacity
market design have primarily been designed with conventional
generation in mind, e.g.\ in GB (\cite{NGCM18}) and in North
America
(\cite{Bowring2013,MoyeMeyn2018}).  Conventional generation is
typically approximated as \emph{firm capacity}, where we define the
latter as idealised capacity which is always available to supply
energy as needed up to a given constant rate.  (To do so nominal
generator capacities are usually multiplied by appropriate
``de-rating'' factors to acknowledge occasional unavailability---see
\cite{NGCM18,Bowring2013}.)

When all capacity is capable of being approximated as if it were firm,
an economic theory of capacity markets is straightforward, and may be
based on balancing procurement cost against cost of unserved energy
(\cite{Stoft2002, Zhao2018}).  \textcite{BothwellHobbs2017} consider
the impact
on societal welfare of the (mathematically) incorrect de-rating of
variable generation,
but consider neither storage nor the mechanism of running a capacity
market.  The need to rapidly reduce fossil fuel dependence means that
both variable generation---e.g.\ wind  and solar power---and
storage now have increasingly important contributions to make to
capacity adequacy (\cite{GeskeGreen}).  The present paper shows how
current approaches to capacity market design
may be extended to give an integrated theory for the inclusion
within a capacity market of all types of capacity provision.  As at
present, the theory is necessarily based on a probabilistic
description of the electricity supply-demand balance process.
However, storage has a natural energy constraint and thus can supply
energy only for a limited period of time before needing to be
replenished; subject to this constraint it may be scheduled flexibly.
Hence, in order to understand both how to schedule storage and to
determine its contribution to capacity adequacy, it is necessary to
pay attention to the sequential statistical structure
of the supply-demand balance process to which that storage is
contributing---see
Sections~\ref{sec:risk-metrics} and~\ref{sec:equiv-firm-capac}.  The
present paper extends and generalises theory which was developed by
the authors in conjunction with National Grid ESO for the integration
of storage into the GB capacity market (\cite{NGECR18}).  However, the
theory is applicable wherever capacity contributions of variable
generation and storage need to be correctly assessed.  This includes
the European and North American markets referenced above.

The determination of a volume of \emph{capacity-to-be-procured} in a
capacity market may be achieved either via the satisfaction of an
appropriate \emph{security-of-supply} standard defined in terms of
some given system \emph{risk metric}, or via the minimisation of an
appropriate economic cost.  (In the latter case the
capacity-to-procure may be variable and specified as a function of the
clearing price in the capacity auction;  this is what currently happens in
GB.)
These two approaches are closely related---see
Section~\ref{sec:economic-approaches}.  In either case, a key step in
the development of an integrated theory
is that of the provision of an appropriate definition of the
\emph{equivalent firm capacity} (EFC) of any capacity-providing
resource.  This EFC is the firm capacity which makes an (appropriately
defined) equivalent contribution to the overall supply-demand balance.
Hence the EFC is necessarily defined with respect to the pre-existing
supply-demand balance process to which this resource is being added
(\cite{ZD2011,DZ2014}).  When the set of capacity-providing resources
contains significant storage, particular care is required in the
determination of EFCs.  One reason for this is the need to account for
the flexibility of scheduling of additional storage added to an
existing supply-demand balance process.  More subtly, when \emph{any}
further resource---including, for example, firm capacity---is added to
an existing set of capacity-providing resources which already contains
storage, that pre-existing storage may \emph{also} be rescheduled so
as to enhance the usefulness of the additional contribution.
A consequence, as we show formally at the end of
Section~\ref{sec:equiv-firm-capac} and demonstrate in the example of
Section~\ref{sec:examples}, is that the EFC of further storage added
to an existing system is \emph{less} (than it would otherwise be) in
the case where that additional system already contains significant
storage.


\label{cr:dr}

Throughout the present paper we treat the process of electricity
\emph{demand} as given.  However, demand response may also be used to
assist in balancing systems.  Demand response has many of the
characteristics of storage, typically making a similarly flexible
contribution.
Its contribution to electricity capacity, and its integration into
capacity markets, may be analysed analogously.  However, in present
day markets it is often treated as demand which may be effectively
foregone---see \cite{LopesAlgarvio} for a review, and also
\cite{NGECR18}.

Sections~\ref{sec:risk-metrics} and \ref{sec:equiv-firm-capac} of the
paper study respectively risk metrics and EFC.  The latter is
necessarily defined in terms of some risk metric and is essential to
the understanding of both capacity adequacy and the operation of a
capacity market.  The studied properties are implicit in the theory of
present markets for what is mostly conventional generation.  However,
so as to understand how to incorporate into such markets both variable
generation and time-limited but flexible resources such as storage, it
is necessary to make these properties explicit.  It is further
necessary to make explicit
assumptions of \emph{continuity} and \emph{smoothness} as available
capacity-providing resource is varied.  We argue that these
assumptions, which are often implicit in other work
(e.g.~\cite{BothwellHobbs2017}), are usually sufficiently satisfied in
practice.  The smoothness assumption yields an important \emph{local
  additivity} property for EFCs; this is essential for the optimal
operation of markets---even in the case where all resource is provided
by firm capacity.  In Section~\ref{sec:equiv-firm-capac} we also show
how to determine the EFCs of marginal contributions of both variable
generation and storage, notably when the objective is the minimisation
of expected energy unserved (EEU).  For storage this requires
consideration of how it may be optimally scheduled.

Section~\ref{sec:capacity-markets} studies the operation of capacity
markets when the objective is that of obtaining at minimum cost
sufficient capacity to meet a given security-of-supply standard
defined in terms of a risk metric.
Section~\ref{sec:economic-approaches} studies the operation of such
markets when the objective is that of the minimisation of an overall
economic cost.  
The present economic theory of such markets requires substantial
modification in the presence of variable generation and, especially,
storage.

The flexibility of storage scheduling has important consequences for
the way in which a capacity market operates, and these are illustrated
in the detailed example of Section~\ref{sec:examples}.  This example
shows the application of nearly all the above theory, and is chosen to
be realistic in the context of a country such as GB.  It further
demonstrates the practical reasonableness of the assumptions required
for a tractable theory.

\vspace{-3ex}

\section{Risk metrics}
\label{sec:risk-metrics}

In the analysis of capacity adequacy, the length of time over which
system risk is assessed---typically a year or a peak season---is
usually divided into~$n$ time periods, each typically of an hour or a
half-hour in length (\cite{BillAll,NGECR18}).  Let random
variables~$D_t$ and $X_t$ denote respectively the total energy demand
and total energy supply in time period~$t$.  Then the
\emph{supply-demand balance} in the time period~$t$ is given by the
random variable $Z_t = X_t - D_t$.  Values of $Z_t$ less than zero
correspond to an energy \emph{shortfall} or \emph{loss-of-load} at
time~$t$.  The \emph{depth} of shortfall at that time is given by
the random variable~$\max(-Z_t,0)$.

Any \emph{risk metric}~$\rho$ is a function of the entire
process $(Z_t,\,t=1,\dots,n)$.  Risk metrics may either be used
directly in the setting of reliability standards---as in the case of
the present GB \lole-based standard
(\cite{NGECR18})---or may arise naturally in \emph{economic}
approaches to determining security-of-supply
(see Section~\ref{sec:economic-approaches}).

\vspace{-2 ex}

\paragraph{LOLE and EEU.}

The two most commonly used risk metrics are \emph{loss-of-load
  expectation} (LOLE) and \emph{expected energy unserved} (EEU), given
respectively by
\begin{align}
  \lole & = \sum_{t=1}^n \Pr(Z_t < 0),\label{eq:1}\\
  \eeu & = \sum_{t=1}^n \E(\max(-Z_t,0))
         = \sum_{t=1}^n \int_{-\infty}^0 \Pr(Z_t < z)\,dz, \label{eq:2}
\end{align}
where~$\Pr$ denotes probability and $\E$ denotes expectation
(\cite{BillAll,KMDHADHSSO,ZD2011}).  It follows from the additivity of
expectations that, in discrete time, LOLE is
the expected number of periods of shortfall during the season under
study, while EEU is the expectation of the sum of the depths of
shortfall during such periods, i.e.\ the expectation of the total
unserved energy.

\label{cr:eeu}
The use of EEU as a measure of economic cost corresponds to a uniform
valuation of unserved energy, regardless of the overall depth of
energy shortfall at any given time and also of the overall duration of
the energy shortfall periods.  The present paper mainly considers such
a uniform valuation.  However, often modest depths or durations of
shortfall may be managed without significant ill effects, while the
avoidance of economic damage becomes increasingly difficult as depth
or duration of shortfall increases.  Thus it may sometimes be natural
to value unserved energy more highly at higher levels of shortfall.

\begin{figure}[ht!]
  \centering
  \includegraphics[scale=0.44]{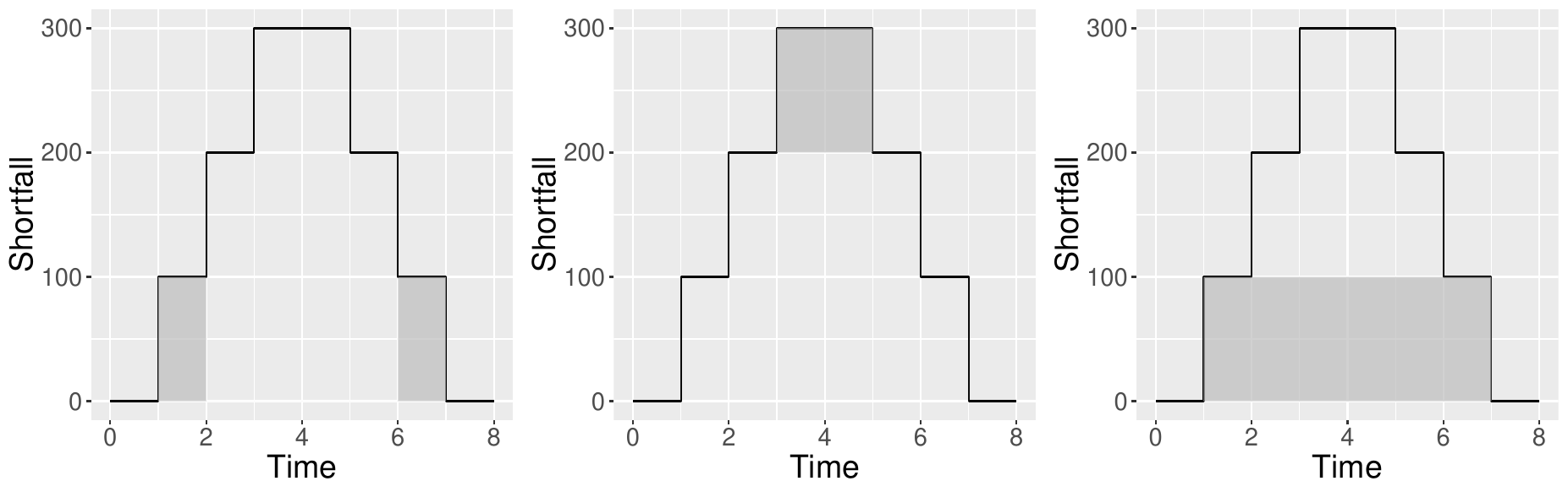}
  \vspace{-1ex}
  \caption{Use of a given volume of storage to mitigate a given period
    of shortfall, and comparison with firm capacity.}
  \label{fig:StorePlots}
\end{figure}

\label{cr:lole}
However, security-of-supply standards are presently more commonly
defined via the use of the LOLE risk metric.  This metric might be
regarded as naturally somewhat flawed, in that---unlike EEU, which
values periods of energy shortfall in proportion to their depth---LOLE
takes no account of shortfall depth.
\label{cr:lole}
Furthermore, we argue below that, in the presence of either variable
generation or storage, LOLE is a poor measure of system reliability.
This is especially so in the case of storage.  We give the relevant
mathematics in Section~\ref{sec:economic-approaches}, but the
following simple example illustrates why the use of LOLE may be
problematic.  Suppose that a certain volume of stored energy is
available to mitigate, but not wholly eliminate, a given period of
energy shortfall.
The loss-of-load duration, and hence LOLE, is minimised by using this
additional stored energy to eliminate entirely the shortfall at those
times within the above period at which this is most efficiently
achieved, i.e.\ at those times at which the \emph{depth} of the
original shortfall is least.  The reason is that this policy maximises
the total length of time for which there is no longer any
loss-of-load---see the left panel of Figure~\ref{fig:StorePlots},
which shows an illustrative pattern of shortfall as above, and the
reduction (shaded) in that shortfall when the stored energy is used to
minimise loss-of-load.  However, it might reasonably be argued that in
practice the available stored energy would be at least as well
employed in instead minimising the \emph{maximum} depth of shortfall
in supply during the period concerned.  The latter policy which would
achieve the same reduction in unserved energy but quite possibly
result in no reduction in loss-of-load whatsoever---see the middle
panel of Figure~\ref{fig:StorePlots} which shows the reduction in
shortfall (again shaded) when the same stored energy is now used in
this way.  Indeed, from an economic perspective, the latter policy
would clearly be better
if the unit cost of unserved energy were an increasing function of
shortfall depth.  Further, regardless of how unserved energy is
valued, the former policy of using stored energy to minimise
loss-of-load
does not make as effective a contribution to system reliability as
would the use of firm capacity to achieve the same reduction in
loss-of-load: unlike storage such firm capacity would continue to
contribute even at those times when shortfall was not completely
eliminated---see the right panel of Figure~\ref{fig:StorePlots} which
shows (again for the same pattern of initial shortfall) the shaded
reduction in shortfall which occurs when firm capacity is used to
achieve the same loss-of-load duration as in the left panel.

The above difficulties in \emph{measuring} the contribution of storage
arise on account of its time-limited duration combined with its
flexibility of use.  In the presence of storage, LOLE may be varied in
a manner which (as illustrated above) does not obviously relate to
system adequacy as usually understood.  To a lesser extent there may
be analogous difficulties with variable generation
(\cite{BothwellHobbs2017});
however, unlike storage, the times at which variable generation is
available are not controllable.  When capacity is provided solely via
conventional generation, and when the latter is capable of being
modelled as if it were (de-rated) firm capacity, then \emph{within a
  given system} the relation of LOLE to an economic measure of system
adequacy such as EEU is essentially predetermined, and the above
difficulties do not then arise.  We give a more precise mathematical
treatment in Sections~\ref{sec:capacity-markets}
and~\ref{sec:economic-approaches}.




\vspace{-2 ex}

\paragraph{General risk metrics.}

We \emph{take as given} the demand process $(D_t,\,t=1,\dots,n)$.
The energy supply process $(X_t,\,t=1,\dots,n)$, and hence the value
of any risk metric, is determined by the set of capacity-providing
resources (conventional generation, variable generation, and storage).
We denote by~$\R$ the set of such resources, and regard any risk
metric~$\rho$ as a function~$\rho{(R)}$ of the set~$\R$.  For a
coherent theory, we assume that the resources within~$\R$ are
optimally used for the minimisation of overall system risk.  This
assumption is needed in the context of flexible resources such as
storage which may be used in different ways to support the
system---however, see Section~\ref{sec:equiv-firm-capac}.  We assume,
as is usual, that any risk metric~$\rho$ is such that $\rho(R)$ is
\emph{decreasing} as the set of resources~$\R$ is increased.

Typically, for a large system, the effect on its overall state---as
measured by the risk metric---of the addition or subtraction of a
single unit of resource is small.  The optimality of the state of the
system (by some criterion) is typically characterised by some form of
\emph{equilibrium} with respect to \emph{marginal}, i.e.\
\emph{relatively small}, variations in that overall state caused by
the addition or subtraction of individual resources.  So as to obtain
a reasonable and tractable economic theory, we require continuity and
smoothness assumptions with respect to such marginal variations.  The
\emph{continuity assumption} is that of reasonably continuous
availability of capacity, i.e.\ that the state of the macroscopic
system
may be varied continuously---at least to a good approximation---by the
addition or subtraction of individual units of capacity-providing
resource.  In large systems this is a reasonable approximation which
is usually made in the design of current capacity markets.  An
exception occurs when considering some very large individual resource,
e.g.\ a large nuclear plant; however, it is straightforward to deal
with a small number of such large contributions---see
Section~\ref{sec:capacity-markets}.  The \emph{smoothness assumption}
is that, for any set of resources~$\R$,
the reduction in risk $\rho(\R) - \rho(\Ri)$ resulting from the
addition of some further marginal resource~$i$ to the set~$\R$ (where,
for simplicity, we denote the resulting set of resources by $\Ri$) is,
to a good approximation, unchanged by \emph{small} variations in the
set~$\R$.  More formally, this is the requirement that, for any small
resource~$j$ (representing the above variation in the set~$\R$) and
for any further small resource~$i$, to a good approximation
\begin{equation}
  \label{eq:3}
  \rho(\Rij) - \rho(\Rj) = \rho(\Ri)  - \rho(\R),
\end{equation}
where, as above, $\Rij$ denotes the set of resources~$\R$ supplemented
by the further resources~$i$ and~$j$, and where, as the contribution
of the resource~$j$ reduces in size, the percentage error in the
relation~\eqref{eq:3} becomes negligible.
This smoothness assumption is essentially a form of differentiability
assumption---see below---and is generally well satisfied in most
applications and for most risk metrics.%
\footnote{The smoothness assumption is not guaranteed: it is possible
  to imagine
  that two capacity-providing resources~$i$ and~$j$ might each make
  identical reductions in risk, as measured by~$\rho$, but be such
  that the use of both together achieved no further reduction in risk
  than the use of either singly, in which case~\eqref{eq:3} would
  fail.} %
In the extended example of Section~\ref{sec:examples}, which is chosen
to be reasonably representative of a system such as that in GB, we
check the applicability of both the above assumptions.

The concept of \emph{firm capacity}, defined as energy supply which is
guaranteed to be available up to a given constant rate throughout the
overall period under consideration, plays a role as a \emph{reference
  measure} in assessing the usefulness of other forms of
resource---see Section~\ref{sec:equiv-firm-capac}.  In considering
variations about the set of resources~$\R$ given by variations in firm
capacity, we also write (analogously to the above) $\R+y$ for the set
of resources~$\R$ supplemented by firm capacity able to supply energy
at a further constant rate~$y$; we further allow that $y$ may be
continuously varied.  This will be important when we consider
\emph{equivalent firm capacity} below.  It is straightforward to show
that, for such variation, the above smoothness assumption implies that
there exists the \emph{derivative} $\rho'(\R)$ of $\rho(\R)$ with
respect to firm capacity; this is such that, for small variations of
the total resource~$\R$ by firm capacity~$y$,
\begin{equation}
  \label{eq:4}
  \rho(\R + y) = \rho(\R) + \rho'(\R) y,
\end{equation}
(the relative error in this approximation tending to zero as $y$ tends
to zero).  This derivative plays an important role in subsequent
analysis.


\vspace{-2ex}

\section{Equivalent Firm Capacity}
\label{sec:equiv-firm-capac}

We take as given a suitable risk metric~$\rho$.  Then, given also any
set of capacity-providing resources~$\R$, the contribution of any
further resource~$i$ (conventional generation, variable generation, or
storage) to be added to the set~$\R$ may be measured by its
\emph{equivalent firm capacity} (EFC) $\efc{i}{\R}$.  This is the firm
capacity which, if added to the set~$\R$ in place of the additional
resource~$i$, would make the same contribution to security-of-supply,
as measured by the risk metric~$\rho$.  Formally, the
constant~$\efc{i}{\R}$ is the solution of
\begin{equation}
  \label{eq:5}
  \rho(\Ri) = \rho(\R + \efc{i}{\R}),
\end{equation}
where, as defined earlier, the notations $\Ri$ and $\R + \efc{i}{\R}$
correspond to the set of resources~$\R$ supplemented respectively by
the resource~$i$ and by firm capacity $\efc{i}{\R}$ (see,
e.g.,~\cite{KMDHADHSSO,ZD2011,NGECR18}).
When the set~$\R$ contains resources such as storage which may be
flexibly used, the addition of further resource may result in a
different pattern of usage of the resources within~$\R$ itself: on
both the left and the right side of~\eqref{eq:5} it is assumed that
the total available resource is being optimally used---see also below
for a discussion of the practicality of this.
When the resource~$i$ is itself firm capacity~$y_i$ for some $y_i>0$,
then $\efc{i}{\R} = y_i$ independently of~$\R$.  However, in general
the EFC of a resource, e.g., variable generation or storage,
\emph{depends also on the existing set of resources~$\R$ to which it
  is being added} (\cite{ZD2011,DZ2014}).  This will be important in
Sections~\ref{sec:capacity-markets} and~\ref{sec:economic-approaches}
for developing satisfactory theories of capacity adequacy and capacity
markets.%
\footnote{It is also possible to define the \emph{equivalent load
    carrying capacity} (ELCC) of any further resource~$i$ with respect
  to an existing set of resources~$\R$.  This is the constant
  $\elcc{i}{\R}$ given by the solution of
  $\rho(\Ri - \elcc{i}{\R}) = \rho(\R)$.  By writing
  $\rho(\R)$ as $\rho(\R - \elcc{i}{\R} + \elcc{i}{\R})$ we have that
  $\elcc{i}{\R} = \efc{i}{\R-\elcc{i}{\R}}$, so that we shall find it
  sufficient in the present paper to work with EFCs.}

It follows from~\eqref{eq:4} and~\eqref{eq:5} that, for the addition
of any \emph{marginal} (i.e.\ small) resource~$i$ to the set~$\R$,
\begin{equation}
  \label{eq:6}
  \rho(\Ri) = \rho(\R) + \rho'(\R) \efc{i}{\R},
\end{equation}
(where the derivative $\rho'(\R)$ in~\eqref{eq:6} remains, as
in~\eqref{eq:4}, that defined with respect to variation in firm
capacity).  Thus, for any existing set of resources~$\R$, and with
respect to the risk metric~$\rho$, the contribution to risk reduction
given by any further marginal resource~$i$ is proportional to its
EFC~$\efc{i}{\R}$.  It hence makes sense to value such additional
resources proportionally to their EFCs---see
Section~\ref{sec:capacity-markets}.  The smoothness
assumption~\eqref{eq:3} further implies the \emph{local additivity} of
the EFCs of marginal (i.e.\ \emph{small}) variations about the
set~$\R$, i.e., for marginal additions $i$ and $j$ to~$\R$,
\begin{equation}
  \label{eq:7}
  \efc{i,j}{\R} = \efc{i}{\R} + \efc{j}{\R},
\end{equation}
where by the left side of~\eqref{eq:7} is meant the EFC of the
combined resources~$i$ and $j$, and where the relative error
in~\eqref{eq:7} again becomes negligible as the additional
resources~$i$ and $j$ reduce in size.  The property~\eqref{eq:7} is a
straightforward consequence of~\eqref{eq:3} and~\eqref{eq:6}:
the substitution of each term in~\eqref{eq:3} by the expression given
by~\eqref{eq:6} yields~\eqref{eq:7} immediately.
\label{cr:la}
The \emph{local additivity} result~\eqref{eq:7} is, like the
smoothness assumption itself, a form of differentiability.  It is a
weaker property than full additivity, i.e.\ it does not hold in the
case of larger variations about the set~$\R$.  This is something which
is well-known and which we discuss further in
Section~\ref{sec:capacity-markets} and illustrate in
Section~\ref{sec:examples}.


\vspace{-2ex}

\paragraph{Determination of EFCs.}
\label{sec:determ-equiv-firm}

The EFC $\efc{i}{\R}$ of any capacity-providing resource~$i$ added to
an existing set of resources~$\R$ is defined by the solution of
equation~\eqref{eq:5}.  The solution of~\eqref{eq:5} may involve trial
values of~$\efc{i}{\R}$ and may not always be straightforward.
Provided that the additional resource~$i$ is marginal (i.e.\ small) in
relation to the existing set~$\R$, then it is generally more
straightforward to obtain the EFC $\efc{i}{\R}$ of the further
resource~$i$ via the solution of~\eqref{eq:6} above, i.e.\ as
\begin{equation}
  \label{eq:8}
  \efc{i}{\R} = \frac{\rho(\Ri) - \rho(\R)}{\rho'(\R)};
\end{equation}
typically all the quantities on the right side of~\eqref{eq:8} may be
readily estimated, e.g.\ via simulation.
(Recall that~\eqref{eq:8} is ultimately a consequence of the earlier
smoothness assumption, the validity of which we check in our extended
example of Section~\ref{sec:examples}.)

As we argue in Sections~\ref{sec:risk-metrics}
and~\ref{sec:economic-approaches}, it is often natural to take the
underlying risk metric~$\rho$ to be given by EEU.  Then, in order to
use~\eqref{eq:8}, we require an expression for the derivative
$\eeu'(\R)$ of~$\eeu(\R)$ with respect to firm capacity.  For any set
of capacity-providing resources~$\R$ which consists entirely of
generation, conventional or variable, we have
\begin{equation}
  \label{eq:9}
  \eeu'(\R) = -\lole(\R).
\end{equation}
The result~\eqref{eq:9} is known in the case where all the resource
within the set~$\R$ is provided by firm capacity
(see~\cite{DECCAnnexC}); its proof in the present case---where~$\R$
may also contain variable generation (but not storage)---is
essentially the same and is given in the appendix.  The
results~\eqref{eq:8} and~\eqref{eq:9} therefore provide an efficient
way to determine the EFCs (with respect to the risk metric given by
EEU) of marginal contributions to capacity in an environment in which
all capacity is provided by generation---a result which is also
implicit in~\textcite{BothwellHobbs2017}.  The result~\eqref{eq:9} is also
important in an economic theory of capacity markets---see
Section~\ref{sec:economic-approaches}.

When the set of resources~$\R$ includes storage, the
result~\eqref{eq:9} requires modification.  The reason for this is
that when further resource---including firm capacity---is added to
an existing set of resources, the use of any storage within that set
\emph{may be rescheduled} so as to continue to contribute optimally to
the minimisation of EEU, thereby increasing the contribution of the
additional resource.  Thus the right side of~\eqref{eq:9} needs to be
replaced by something which is larger in absolute value.

\label{cr:simpopt}
For an exact analysis, assume that storage may be completely recharged
between periods of what would, in the absence of storage, be energy
\emph{shortfall}, but that storage may not usefully be recharged
during periods of shortfall.  (This is currently the case in GB, for
example, where periods of shortfall are generally well separated;
typically there is at most a single period of shortfall during any
day, and storage may be fully recharged overnight.)  Assume also that
any process of what would otherwise be continuous energy shortfall is
met as far as possible from available storage, with the aim of
minimising residual unserved energy, and that each store~$i$ delivers
energy subject to \emph{rate} (power) and total \emph{energy}
constraints.%
\footnote{In other contexts the energy constraint might be referred to
  as the ``capacity'' of the store (as measured, for example, in MWh).
  However, in the context of capacity markets, ``capacity'' has units
  of power (as measured, for example, in MW), and we preserve this
  usage throughout this paper.} %
\label{cr:salg}
At each time~$t$ define the \emph{residual lifetime}
of each store~$i$ as such residual energy in the store at time~$t$ as
is available for subsequent use
divided by the maximum rate at which that energy can be served.  Then
the above minimisation of EEU is achieved by the use of the
\emph{greedy} algorithm in which, at each successive time~$t$ during
any period of shortfall, this shortfall is reduced as far as possible
from energy in storage and in which the use of the stores is
prioritised in descending order of their residual lifetimes.
The optimality of the above policy is proved by~\cite{CZ2018} and
by~\cite{ETA-pscc}.
\label{cr:fk}
Note that the implementation of this policy, although optimal,
nevertheless does not, as each successive time~$t$, require any
foreknowledge of the shortfall process subsequent to time~$t$, so that
this policy is practicable to the extent that the scheduling of stores
may be coordinated.  (If, in more extreme situations, storage may not
be fully recharged between periods of use, then optimal scheduling is
more complicated and does require foresight as above.  This is also
the case when the unit valuation of unserved energy is not uniform but
increases according to the depth of shortfall, in which case it is
sometimes optimal to defer serving energy from storage in anticipation
of its being more usefully supplied at a later
time---see~\cite{CZ2018} for a detailed analysis.)
\label{cr:algs}
Even when such coordination is difficult, under the above assumption
that storage may not be recharged within shortfall periods, any policy
for the use of the stores which does not actually spill energy will
clearly also perform optimally in any shortfall period in which either
the stores are emptied or the shortfall eliminated, so that most
reasonable policies for the use of storage may generally be expected
to work well---again see \cite{CZ2018}.
\label{cr:sedef}
Let~$\S_e$ be the set of
stores which, under the above optimal policy, are empty (i.e.\ have
reached their minimum permitted state of charge) at the end of the
shortfall period.  Then it is formally shown by~\cite{CZ2018}
that, in the presence of storage, the result~\eqref{eq:9} is replaced
by the more general result
\begin{equation}
  \label{eq:10}
      \eeu'(\R) = -\lole(\R\setminus\S_e),
\end{equation}
where $\R\setminus\S_e$ is the set of resources in~$\R$ other than
those in the set~$\S_e$.
\label{cr:simul}
(Since the shortfall process is typically random, so also is the set
$\S_e$.  Thus if, for example, the right side of~\eqref{eq:10} were
being evaluated by simulation, for each such simulation the set $\S_e$
and the loss-of-load corresponding to $\R\setminus\S_e$ would both be
observed; the latter would then be averaged over simulations to
estimate $\lole(\R\setminus\S_e)$.)



\label{cr:toy}
As observed above, when firm capacity is added to a set of
resources~$\R$, it contributes more effectively to the reduction of
EEU when the set~$\R$ already contains storage resources which may
then be rescheduled.  As a simple example, consider a shortfall of
100~MW for one hour followed by 200~MW for a further hour.  A store of
100~MWh capacity with a rate of 100~MW can contribute its 100~MWh
during the first hour, thereby eliminating the shortfall during that
hour.  Firm capacity of 100~MW can now contribute to a further
reduction in unserved energy of 100~MWh, achieved during the second
hour.  However, if the earlier storage is rescheduled to contribute
(equally effectively) during the second hour, the subsequent addition
of the 100 MW of firm capacity eliminates the shortfall entirely and
achieves a further reduction in unserved energy of 200 MWh (thereby
achieving a bonus of 100~MWh from the storage rescheduling).

A corollary is that the EFC of further \emph{storage resources} is
typically \emph{reduced} by the presence of existing storage resources
in a system---see Section~\ref{sec:examples}.

\vspace{-2ex}

\section{Capacity markets}
\label{sec:capacity-markets}

In this section we consider the operation of a capacity market where
the objective is that of obtaining at minimum cost a sufficient set of
capacity-providing resources~$\R$ so as to satisfy some condition
\begin{equation}
  \label{eq:11}
  \rho(\R) \le k,
\end{equation}
expressed in terms of some appropriate risk metric~$\rho$.  We assume
throughout that there is a fixed process of net demand to be met by
these resources.  The level~$k$ may either be chosen so that the
condition~\eqref{eq:11} defines some given security-of-supply
standard, or it may be chosen according to some economic criterion
(see Section~\ref{sec:economic-approaches}).
We allow that some resource may be provided by facilities other than
firm capacity, for example variable generation or storage.  The
required theory is the same as that which is currently used when all
resource (typically conventional generation) may be approximated as
firm capacity, except that it is now necessary to appropriately define
the EFCs of other resources.  Further some additional calculation may
be required in the operation of any auction associated with the
capacity market---at least when resource other than firm capacity
makes a substantial contribution.

\label{cr:gbcm}
An example is given by the capacity market which currently operates in
GB (\cite{NGCM18}), and which typically seeks to secure at minimum
cost a level of capacity compatible with the GB LOLE-based
security-of-supply standard (however, see also
Section~\ref{sec:economic-approaches}).  The auction associated with
this market determines a unit \emph{clearing price} such that any
successful bidder (capacity provider) receives this clearing price
multiplied by its offered capacity.  When the latter is other than
firm capacity, e.g.\ variable generation or storage, the bidder's EFC
is used instead and, since the total such capacity is currently small,
is reasonably estimated in advance of the auction---essentially as
described in Section~\ref{sec:equiv-firm-capac}.
The GB auction is then conducted as a (single clock) \emph{descending
  clock auction}.  An initially high unit price is gradually reduced;
bidders may offer discrete capacity resources and may exit the auction
at any point; the auction is cleared at the point where there is
just sufficient capacity remaining in the auction to meet the required
reliability standard~\eqref{eq:11} for the given LOLE-based risk
metric~$\rho$ (a process which clearly involves running the auction to
just beyond this point); the unit price at this point then becomes the
clearing price, and each bidder remaining in the auction is then paid
as above.  For more details, see \cite{NGCM16,NGCM18}.  The auctions
associated with other systems, e.g.\ in North America, may be more
complex as these systems may be partially fragmented by capacity
constraints between different regions---see, e.g., \cite{ISONE}---and
require multiple clocks, but the general theory below remains
applicable.

Such a descending clock auction relies on the EFC of any capacity
other than firm capacity being clearly identified in advance.
However, as discussed in Section~\ref{sec:equiv-firm-capac}, the EFC
of any such capacity-providing resource depends also on the overall
supply-demand balance process of which it forms a part---which is that
determined by the overall set of resources~$\R$ finally selected in
the auction.
It follows that, in the presence of a \emph{substantial} number of
resources other than firm capacity, the form of auction described
above may require some adjustment as described below.

In order to provide a better basis for a coherent theory, instead of
considering the minimum \emph{unit} price (i.e.\ price per unit of
EFC) which each bidder is prepared to accept for its capacity
offering, we consider instead the minimum \emph{total} price~$c_i$
which each bidder or resource~$i$ is prepared to accept in return for
making available some given capacity.  (It is natural that bidders
should have such a price in mind, and, in the theory below there is
then no need for bidders to estimate themselves the EFCs of the
resources they are offering.)  For example, this capacity might be a
given level of firm capacity for as long as might be required, or it
might be storage which could be called upon as required and which
could be used flexibly subject to specified rate and energy
constraints.  The societal problem is now to design an auction to
choose at minimum cost a subset $\R$ of those resources competing in
the market, such that the required reliability condition~\eqref{eq:11}
is satisfied.  We continue to assume that such an auction is
pay-as-clear, i.e.\ associated with its outcome is a unit clearing
price~$p$ such that, if~$\R$ is the set of resources which are finally
successful in the auction, then
\begin{equation}
  \label{eq:12}
  \begin{split}
    c_i \le p \times \efc{i}{\R}, & \qquad i \in \R,\\
    c_i > p \times \efc{i}{\R}, & \qquad i \notin \R,
  \end{split}
\end{equation}
and each successful resource~$i\in\R$ is paid in total
$p \times \efc{i}{\R}$.  The relations~\eqref{eq:12} define a
competitive equilibrium condition necessarily satisfied by the unit
clearing price~$p$ and the required optimal set of resources~$\R$.
That this is so depends implicitly on the continuity and smoothness
assumptions of Section~\ref{sec:risk-metrics}.  The continuity
assumption ensures that resource is reasonably continuously variable,
and the smoothness assumption guarantees the local additivity property
of Section~\ref{sec:equiv-firm-capac}; under these assumptions, were
such a unit clearing price \emph{not} to exist for the claimed optimal
set of resources~$\R$, then resource from outside that set might be
more cheaply substituted for resource from within it while continuing
to satisfy the required reliability condition~\eqref{eq:11},
contradicting the optimality of the set~$\R$.  In general we might
then expect the relations~\eqref{eq:11} and~\eqref{eq:12} to define
the required set of resources~$\R$ uniquely:
\label{cr:converge}
if the EFCs of the
individual resources were fixed and if the resource provided by each
bidder were continuously variable, then it is clear that the auction
clearing mechanism would indeed identify the minimum-cost set of
resources satisfying the constraint~\eqref{eq:11}; for $\R$ in the
neighbourhood of this minimum, the EFCs $\efc{i}{\R}$ vary only slowly
with~$\R$ so that the same result might reasonably be expected to
hold; however, the fact that in reality resource is offered to the
auction in discrete quantities means that, because of minor
``overshoot'' problems, the minimum-cost resource set may not be
precisely identified.  When, in the presence of one or more
\emph{large} capacity-providing resources, the above continuity
assumption breaks down, there may be more major problems of
``overshoot''.  However, these problems may usually be solved with a
little experimentation; otherwise an integer optimisation approach is
required.  See also the example of Section~\ref{sec:examples} here.

 

When the EFCs $\efc{i}{\R}$ may be reasonably be estimated in advance,
the unit clearing price~$p$ may be identified by a
descending clock auction as described above.
Alternatively, if bidders were willing to declare their minimum total
prices~$c_i$ in advance, an auction could be conducted offline by
ranking in ascending order the minimum unit
prices~$c_i/\efc{i}{\R}$---to give what is usually referred to in the
context of energy markets as a \emph{merit order stack}
(\cite{Ofgem2016,StaffellGreen2016}); the unit clearing price~$p$
would then be chosen so that the accepted set of resources~$\R$
satisfying~\eqref{eq:12} was just sufficient to satisfy the required
reliability condition~\eqref{eq:11}.

\label{cr:auctioneer}
However, when there are substantial resources other than firm capacity
participating in the market, the final EFCs~$\efc{i}{\R}$ of these
resources (estimated with respect to the finally accepted set~$\R$)
may \emph{not} be sufficiently known in advance of any capacity
auction.  In order to identify a unit clearing price~$p$ and resource
set~$\R$ such that the required conditions~\eqref{eq:11}
and~\eqref{eq:12} hold, it may be necessary to proceed iteratively:
starting with initial estimates (by the auctioneer) of the unknown
EFCs, an initial clearing price~$p$ and initial accepted resource
set~$\R$ may be obtained---e.g.\ through the formation of a
merit-order stack as above; on the basis of this set~$\R$, EFCs may be
re-estimated and improved values of the clearing price~$p$ and
accepted resource set~$\R$ may then be obtained as before;
one might reasonably expect convergence to a solution of~\eqref{eq:11}
and~\eqref{eq:12} within a fairly small number of iterations---again
see the example of Section~\ref{sec:examples}.

\label{cr:lm}
Note that, in the above theory, once the EFCs of the various offered
resources are settled, the auction clearing price has an
interpretation as a shadow price (or Lagrange multiplier) associated
with the requirement constraint~\eqref{eq:11} of the minimum-cost
optimisation problem.  This interpretation opens the way to the use of
optimisation theory more generally to solve analogous problems with
multiple constraints---as in the case of auctions for multiple
products.


Finally, recall also that, while the theory of this section depends on
the local additivity of EFCs given by~\eqref{eq:7}, such additivity
does not hold over more extensive variation of a set of
resources---again see Section~\ref{sec:examples}.

\vspace{-1ex}

\section{Economic approaches}
\label{sec:economic-approaches}

Section~\ref{sec:capacity-markets} considered the determination of the
optimal set~$\R$ of capacity-providing resources meeting a given
security-of-supply standard.  However, it is also possible to consider
explicitly economic approaches to the determination of electricity
capacity.  Thus one might choose the set of resources~$\R$ so as to
minimise an overall economic cost
\begin{equation}
  \label{eq:13}
    \voll\times\eeu(\R) + c(\R),
\end{equation}
where the constant $\voll$ is a unit \emph{value-of-lost-load},
$\eeu(\R)$ is, as before, the expected energy unserved associated with
the optimal use of the set of resources~$\R$, and $c(\R)$ is the cost,
within a capacity market, of providing the set of resources~$\R$.
This measure is common in economic approaches to the determination of
electricity capacity adequacy
(\cite{BillAll,BothwellHobbs2017,DECCAnnexC}).
We examine the extent to which this approach generalises to
include variable generation and storage.  Of particular interest is
the extent to which an economic criterion based on the valuation of
EEU may continue to be reduced to a risk-based criterion expressible
in terms of an LOLE-based security-of-supply standard.  It turns out
that this may or may not be possible in the case where the
capacity-providing resources include variable generation, depending on
the statistical characteristics of the latter, but is not readily
possible when these resources include storage.  We thus consider three
cases of increasing generality.

When all resource may be approximated as firm capacity, and when this
is reasonably continuously available, the set of resources~$\R$ may be
identified with the level of firm capacity it provides.  The overall
cost~\eqref{eq:13} is a convex function of that capacity, and is
minimised at the level of capacity~$\R$ such that
\begin{equation}
  \label{eq:14}
  \voll\times\eeu'(\R) + c'(\R) = 0,
\end{equation}
where, as usual, $\eeu'(\R)$ is the derivative of $\eeu(\R)$ with
respect to firm capacity, evaluated at~$\R$, and where $c'(\R)$ is the
similarly the derivative with respect to firm capacity of the
cost~$c(\R)$ of obtaining that capacity.  In GB the
derivative~$c'(\R)$, evaluated at the optimal level of capacity~$\R$
as given by the solution of~\eqref{eq:14}, is commonly referred to as
the \emph{cost of new entry} (CONE) (\cite{DECCAnnexC}).  It thus
follows from~\eqref{eq:14} and from the earlier result~\eqref{eq:9}
that the optimal level of capacity minimising~\eqref{eq:13} is given
by the solution of
\begin{equation}
  \label{eq:15}
  \lole(\R) = \frac{\cone}{\voll}.
\end{equation}
The quantity~$\cone=c'(\R)$, evaluated at the optimal level of
capacity~$\R$, is of course determined by the capacity market.
However, $c'(\R)$ usually varies slowly with~$\R$ and may often be
estimated in advance of any capacity auction---e.g.\ on the basis of
earlier auctions (\cite{DECCAnnexC}).
The relation~\eqref{eq:15} then re-expresses the economic criterion of
the present section (that of minimising~\eqref{eq:13}) as a simple
LOLE-based criterion, and the determination of the minimum-cost level
of capacity~$\R$ such that~\eqref{eq:15} holds is as described in
Section~\ref{sec:capacity-markets}.

The above theory is sometimes used as an economic justification for
the present LOLE-based GB reliability standard.  A description is
given in~\textcite{DECCAnnexC}, where the central values of $\voll$
(£17/kWh) and $\cone$ (£49/kW per year) identified there correspond
approximately, via~\eqref{eq:15}, to the present GB reliability
standard of a (maximum) LOLE of 3 hours per year.  Nevertheless, the
values of $\cone$ identified in recent GB capacity auctions have
varied widely and are in all cases less than the value quoted above
(see~\cite{NGCM18}); thus the above economic justification for the
present GB standard remains a somewhat theoretical one.  An
alternative mechanism, now implemented within the most recent GB
capacity auctions, is that of including provision for obtaining a
total level of capacity which depends on the bidding taking place
within the auction itself, so that the lower the clearing price the
greater the capacity obtained---in greater conformity with the
solution of~\eqref{eq:15}.  This may be done through the specification
of a \emph{demand curve}
(\cite{NGEMRDC,MoyeMeyn2018,Zhao2018,Hobbsetal2007}) which specifies
the total level of capacity to be obtained as a function of the
clearing price in the auction.  This approach is straightforwardly
implemented in, e.g., the GB descending clock auction as described in
Section~\ref{sec:capacity-markets}: as the unit price is decreased in
successive rounds of that auction, so the target capacity to be
obtained is increased in line with the specified demand curve; the
auction clears when the offered capacity remaining in the auction is
equal to the current target capacity.  While the GB demand curve is
currently determined by government policy, it could of course be
chosen so that the relation~\eqref{eq:15} was satisfied for the set of
resources~$\R$ obtained in the capacity auction, with CONE as the
dynamically determined clearing price of the auction.  We note also
that the GB capacity market does not attempt to take account of
possible variation in \emph{energy} costs according to the type of
capacity selected; in practice this would be extremely difficult.

Of interest now is the extent to which the above theory generalises to
the more complex situation in which all capacity is provided by either
conventional or variable generation, but in which storage is not
present.  There is then no scalar measure of capacity which is
sufficient to determine either EEU or LOLE.  We do, however, have the
following result.

\begin{proposition}\label{prop:one-one}
  Suppose that all capacity is provided by some form of generation,
  and that this is reasonably continuously variable (see
  Section~\ref{sec:risk-metrics}).  Suppose further that, as the
  set~$\R$ of such resources is varied, there is a one-one
  correspondence between values of $\eeu(\R)$ and values of
  $\lole(\R)$.  Then the optimal set of resources~$\R$ minimising the
  overall economic cost~\eqref{eq:13} is again that which minimises
  the cost~$c(\R)$ of providing them subject to the
  constraint~\eqref{eq:15}.
\end{proposition}

A formal proof of Proposition~\ref{prop:one-one} is given in the
appendix.  The essence of the argument is that the existence of the
above one-one correspondence ensures that minimisation of $c(\R)$
subject to a constraint on $\eeu(\R)$ is equivalent to minimisation of
$c(\R)$ subject to the corresponding constraint on $\lole(\R)$.
However, the above correspondence between values of $\eeu(\R)$ and
$\lole(\R)$, while clearly trivial in the case where all resource is
provided by firm capacity only (since both are decreasing functions of
the level of that firm capacity), is not \emph{guaranteed} in the case
where resource is also provided by variable generation.  It is
possible that, in assessing its contribution to capacity adequacy,
variable generation is capable of being treated \emph{as if} it were
firm capacity at a constant and appropriately ``de-rated'' level---as
when the process of variable generation is statistically independent
of that of demand---so that the above one-one correspondence between
EEU and LOLE is maintained (see~\cite{ZD2011}).  However, it is also
possible that the pattern of availability of variable generation
is such that this correspondence is not maintained.
\label{cr:oneone}
For example, it is possible that on occasions periods of solar
generation might be contained within periods of loss-of-load in such a
way that this generation contributes to the reduction of unserved
energy without reducing at all the duration of the loss-of-load
periods.  (This is particularly possible in countries where energy
shortfalls tend to occur in the middle of the day.)  The determination
of the optimal set of capacity-providing resources minimising the
overall economic cost~\eqref{eq:13} may then require numerical
analysis.

Finally, the above theory does not directly generalise to the case
where the capacity provision also includes storage, for then, for
the reasons indicated in Section~\ref{sec:equiv-firm-capac},
the result~\eqref{eq:9}---upon which the above theory depends---no
longer holds.  In this case the determination of the optimal set of
capacity-providing resources minimising~\eqref{eq:13} may again
require numerical analysis.

\vspace{-1ex}

\section{Example}
\label{sec:examples}

We present a detailed example of energy storage and firm capacity
competing in a capacity market and designed to illustrate most of the
theory of the present paper.  The dimensions of the example correspond
approximately to those of the current GB electricity system, except
that we allow more storage than is currently present.  The example
illustrates how to value the contribution of individual stores within
a market to which storage makes a significant contribution, so as to
ensure the optimal operation of this market.  It does not include
variable generation, which does not currently participate in the GB
capacity market; however, we emphasise that variable generation may be
handled and valued in exactly the same way as illustrated here for
storage.  Since substantial storage is involved, we take as our
objective that of choosing at minimum cost a set of resources~$\R$ to
meet a given EEU reliability standard (see
Sections~\ref{sec:risk-metrics} and \ref{sec:economic-approaches}).
We take this to be 2746~MWh per year, as in the current GB system with
relatively little storage this corresponds to an LOLE of approximately
3 hours per year (the current GB standard).

\label{cr:bgd}
We first create a credible \emph{background} supply-demand balance
process against which the capacity auction is to take place---it is
typically the case in GB, for example (but not necessarily elsewhere),
that in any auction there is pre-existing capacity already committed,
e.g.\ from multi-year auctions held in previous years, which therefore
does not compete in the current auction.  For this background process
we assume a set of 230 conventional generators with a total of
61.36~GW of installed capacity.  Capacities and outage probabilities
for these conventional generators correspond approximately to a
2015--16 National Grid scenario for GB.
The availability of each generator is modelled as a two-state Markov
process in which each generator is either completely available or
completely unavailable, with a mean time to repair of 50 hours and a
mean time to failure such that the equilibrium outage probability of
the generator agrees with National Grid's scenario values (see
\cite{ESDT}).
\label{cr:illus}
For the purposes of this illustrative example, an
empirical demand-net-of-wind process is created from paired hourly
observations of GB electricity demand and wind generation for the
winter season 2010--11 rescaled to 2015--16 levels of demand and an
assumed installed total wind generation capacity of 14~GW (see
\cite{WZD2018}).  (In practice a longer demand-net-of-wind series---of
approximately 10 years---is used in the GB capacity assessment
analysis, as this process varies considerably from year to year.)
From this demand-net-of-wind process are subtracted 100 simulations of
the conventional generation process to create 100 simulations of a
residual demand process, which defines a sufficiently representative
background process for the present example.  The further \emph{firm
  capacity} which would require to be subtracted from this background
demand process in order to meet the target reliability standard of
2746~MWh per year is 1973~MW.  However, this residual demand or
background process is to be managed instead from the further
generation and storage resources obtained in the capacity market.  The
volume of resource to be thus obtained corresponds to that which might
be required in a ``top-up'' capacity market, such as that held in GB
one year ahead of real time.

Competition in the capacity market is provided by an assumed set of
120 stores and 30 units of firm capacity.  The rate and energy
constraints of the stores are chosen to be representative of those
typically found in systems such as that of GB.  However, in order to
illustrate some of the concerns most relevant to future systems in
which storage may play a larger part, we choose a relatively large
number of stores.  We assume there are 60 stores with a rate
constraint of 50 MW; 10, 15, 15 and 20 of these stores have energy
constraints of respectively 12.5, 25, 50 and 100 MWh.  We further
assume there are 60 stores with a rate constraint of 100 MW; 10, 15,
15 and 20 of these stores have energy constraints of respectively 25,
50, 100 and 200 MWh.  The firm capacity units are assumed to have
capacities between 10 MW and 100 MW. There are three units for each
multiple of 10 MW capacity (i.e.\ three units with capacity of 10 MW,
three with capacity of 20 MW, etc).

The minimum total prices~$c_i$ at which the stores or unit of firm
capacity~$i$ are prepared to make themselves available (see
Section~\ref{sec:equiv-firm-capac}) correspond, for firm capacity, to
a range of \emph{unit} prices with a mean of \pounds 32/kW and
standard deviation of \pounds 4/kW, while those of storage correspond
to unit prices with a mean of \pounds 30/kW and standard deviation of
\pounds 3/kW.  (The latter storage prices are per unit of EFC
$\efc{i}{\R}$ calculated with respect to the set of resources finally
obtained in the auction---while the values of the $c_i$ for all
resources are, of course, specified at the outset of the auction.)
These choices of the $c_i$ were made to ensure reasonable competition
in the market between firm capacity and storage.

The EFC $\efc{i}{\R}$ of any storage unit~$i$, relative to
any set of resources~$\R$ to which it is being added, is calculated as
in Section~\ref{sec:equiv-firm-capac}.  In all cases storage is
optimally scheduled for the minimisation of EEU as described in that
section on the assumption that all storage may be completely recharged
overnight, but not at other times.

\begin{figure}[!ht]
  \centering
  \includegraphics[scale=0.6]{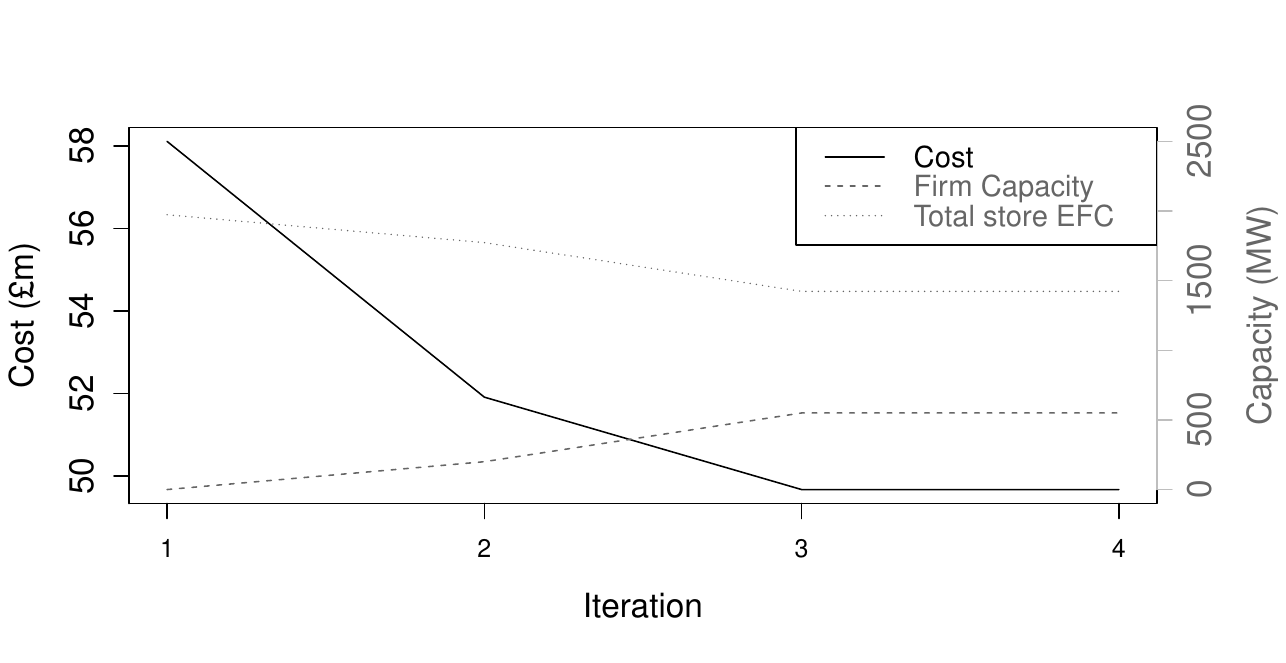}
  \vspace{-3 ex}
  \caption{Iterative procedure: total market cost and associated
    capacity mix at the end of successive iterations.}
  \label{fig:cvge}
\end{figure}

\label{cr:fpt}
The determination of the optimal set of resources~$\R$ meeting the
required EEU reliability standard at minimum cost is as described in
Section~\ref{sec:capacity-markets}.  It is necessary to identify~$\R$
such that the conditions~\eqref{eq:11} and~\eqref{eq:12} are
satisfied.  We take a ``fixed-point'' approach.  Thus it is necessary
to determine a set of resources~$\R'$,
such that, if the EFCs $\efc{i}{\R'}$ of all resources~$i$
participating in the capacity market are calculated with respect
to~$\R'$, and if then a set of resources~$\R$ such that the reliability
condition~\eqref{eq:11} and the equilibrium conditions~\eqref{eq:12}
(with $\efc{i}{\R}$ replaced by $\efc{i}{\R'}$ in~\eqref{eq:12})
are satisfied for some clearing price~$p$, then we have $\R'=\R$.
We therefore proceed iteratively; each iteration starts with an input
mix of resources~$\R'$ satisfying~\eqref{eq:11}, and then an output
mix of resources~$\R$ is determined as above; convergence occurs when
$\R'=R$ and~\eqref{eq:11} and~\eqref{eq:12} are then satisfied for the
final clearing price~$p$.  We take the initial input set~$\R'$ to
consist entirely of firm capacity.  For each subsequent iteration a
simple choice of input set~$\R'$ would be to take this to be the
output set from the previous iteration.  However, the competitiveness
of storage and firm capacity in the present example means that
convergence is oscillatory and relatively slow.  It is considerably
accelerated by instead taking $\R'$ to consist of firm capacity equal
to the mean of the input and output firm capacities associated with
the previous iteration, supplemented by sufficient storage chosen in
accordance with the merit-order stack resulting from the previous
iteration such that the reliability condition~\eqref{eq:11} is
satisfied.  The results are as shown in Table~\ref{tab:results}, which
shows in particular the associated total cost and auction clearing
price at the end of each iteration, while Figure~\ref{fig:cvge} plots
total cost, total firm capacity and sum of storage EFCs at the end of
each successive iteration.


\begin{table}[ht]
  {\small
    \begin{tabular}{|c|cccccc|}
      \hline
      Iteration&Total cost & Clearing&LOLE & EEU&  Sum of store& Firm capacity\\
               & (\pounds m)&price (\pounds/MW) &(hrs/year)& (MWh) & EFCs (MW)&(MW)\\
      \hline
      1&58.1&19,379&3.56&2745.8&2999  &0\\
      2&51.9&27,353&3.26&2745.8&1698  &200\\
      3&49.7&29,492&2.92&2744.5&1134.4&550\\
      4&49.7&29,492&2.92&2744.5&1134.4&550\\
      \hline
    \end{tabular}
  }
  \caption{Iterative procedure: total market cost and associated
    outputs at the end of successive iterations.}
  \label{tab:results}
\end{table}

\vspace{-1ex}

Convergence to a set of resources~$\R$ such that \eqref{eq:11}
and~\eqref{eq:12} are then satisfied for the associated clearing
price~$p$ is obtained after three iterations---as confirmed by a
fourth.  Recall that, subject to the continuity and smoothness
conditions of Section~\ref{sec:risk-metrics}, these conditions are
necessarily satisfied at the optimal set~$\R$, and---as argued in
Section~\ref{sec:capacity-markets}---might reasonably be expected to
identify it uniquely.  However, to the extent that capacity offerings
in the market are discrete, absolute optimality cannot be
guaranteed. (The final set of resources obtained here is identical
with that obtained by the slower, simpler algorithm also discussed
above.)  The solution obtained---the set of resources~$\R$
meeting the required reliability standard at minimum cost on a
\emph{pay-as-clear basis}---is a combination of 550 MW of firm
capacity and a set of stores the sum of whose EFCs $\efc{i}{\R}$
\emph{evaluated with respect to~$\R$} is 1134 MW.  These are the EFCs
appropriate to the \emph{marginal} contributions of the individual
stores at the point where the market clears (and as is appropriate to
the optimal operation of the capacity market---see
Section~\ref{sec:capacity-markets}), and determine their payments
received in the capacity market.
However, the EFC of the entire set of accepted stores---treated as a
single unit---is 1423 MW, i.e.\ this is the amount of further firm
capacity which would be required to substitute for the entire set of
accepted stores in order to meet the required EEU reliability
standard.  (The total power rating of these stores is 3700 MW, but the
store durations are quite low.)  This is an extreme case
of the \emph{nonadditivity} of EFCs over significant numbers of
resources, as discussed at the end of
Section~\ref{sec:capacity-markets}.  It is further a reflection of the
observation, at the end of Section~\ref{sec:equiv-firm-capac}, that
the EFC of further storage resources is typically reduced by the
presence of existing storage resources.  A similar phenomenon is well
known in the case of some variable generation such as wind power,
although the reason for the latter is very different and results
from the temporal coincidence of wind resources at different
locations.

In GB (and other countries) the capacity market is settled through a
single-pass descending-clock auction which identifies the required
clearing price.  Thus the EFCs of storage facilities---which now
participate in the GB capacity market---are estimated in advance of
the capacity auction.  The EFC $\efc{i}{\R}$ of each store~$i$ is
determined as described in the present paper, but---as for the first
iteration of the present example---is done so with respect to an
``initial'' set of resources~$\R$ which is taken to be the amount of
\emph{firm capacity} which would be required in order to meet the GB
reliability standard.  In GB most resource currently participating in
the market is conventional generation, and the storage EFCs estimated
as above are close to their true values, i.e.\ to those calculated
with respect to the set of resources~$\R$ finally obtained in the
market.  However, the example of the present section is chosen so that
storage plays a significant role---something which may very well also
be the case in future real systems.  In this example the initial EFC
estimates $\efc{i}{\R}$ of the stores, determined as described above,
prove to be considerable overestimates in comparison with their true
values.  The reason for this is as discussed at the end of
Section~\ref{sec:equiv-firm-capac}: storage added to existing storage
is less valuable than when added to firm capacity providing the same
level of reliability.  A consequence, \emph{in the present example},
of this overvaluation of the contribution of storage would be that, if
uncorrected, \emph{all} the resource obtained in the capacity market
would consist of storage---as at the end of the first iteration.
Further, with the realistic resource costs chosen for this example,
the cost of obtaining sufficient (all storage) resources to correctly
meet the required reliability standard is \pounds 58.1m, whereas the
true cost of the optimal resource set (as identified earlier and
consisting of a mixture of firm capacity and storage) to meet that
reliability standard is \pounds 49.7m.  The more careful evaluation of
the contribution of storage in the present example therefore leads to
a cost saving of 14.5\%.

\begin{figure}[ht!]
  \centering
  \includegraphics[scale=0.45]{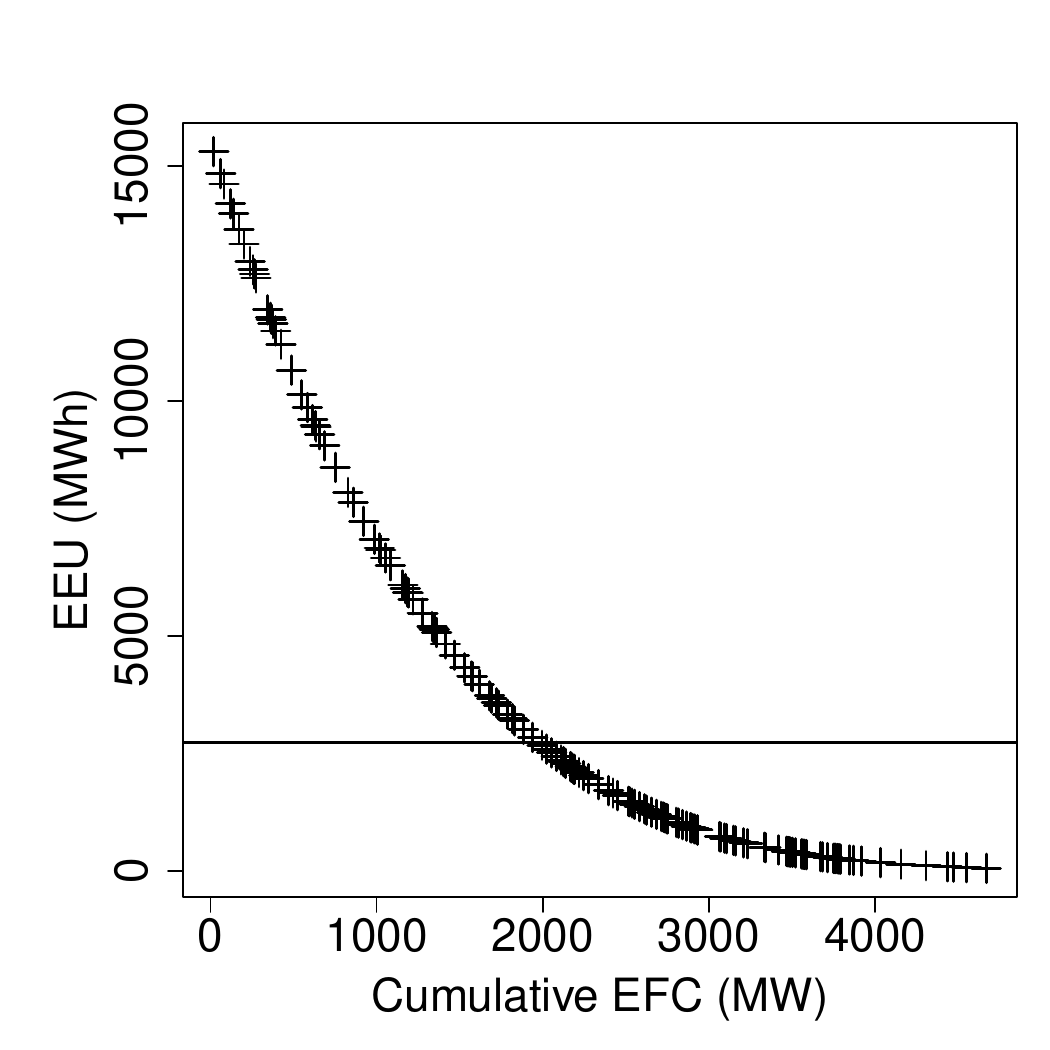}
  \vspace{-1ex}
  \caption{Plot of residual EEU again cumulative EFC.}
  \label{fig:continuity}
\end{figure}

Finally we test more carefully the extent to which the continuity and
smoothness assumptions of Section~\ref{sec:risk-metrics} are
applicable in the current example.  Figure~\ref{fig:continuity} shows
the effect of starting with the background system of pre-existing
capacity and gradually adding to it all the capacity-providing
resources competing in the auction, taken in the order of the final
merit-order stack.  The figure plots residual EEU against cumulative
EFC.  At each point the latter is again the firm capacity which would
substitute for the resources added so far while maintaining the same
level of residual EEU (so that the underlying relationship defined by
the plotted points is in fact independent of the order in which the
resources are taken).  We see that, as required for the continuity
assumption, there are no large gaps between successive points.  In
particular these points are dense in the region corresponding to the
target EEU for the auction---as represented by the horizontal line.
We therefore conclude that the continuity assumption is sufficiently
well satisfied.


\begin{figure}[ht!]
  \centering
  \includegraphics[scale=0.6]{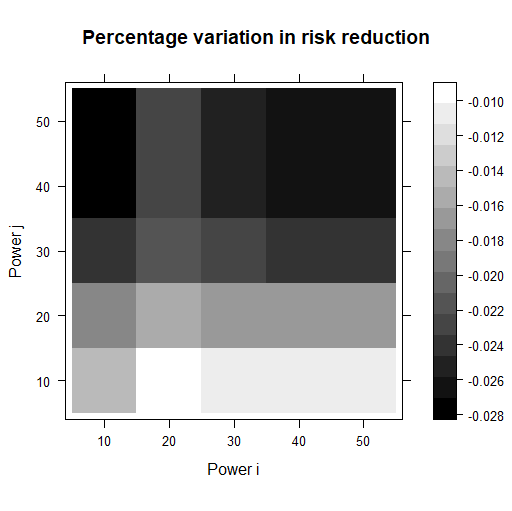}
  \vspace{-4ex}
  \caption{Heat plot to examine the accuracy of the smoothness
    assumption~\eqref{eq:3}.}
  \label{fig:smoothness}
\end{figure}

Figure~\ref{fig:smoothness} examines the validity of the smoothness
assumption~\eqref{eq:3}.  It shows, as a ``heat plot'', the percentage
variation
\label{cr:monotonicity}
\begin{displaymath}
       \frac{\rho(\Rij) - \rho(\Rj)
          - (\rho(\Ri)  - \rho(\R))}%
        {\rho(\Ri)  - \rho(\R)}
        \times 100\%
\end{displaymath}
between the left and right sides of~\eqref{eq:3}, where the
set~$\R$ is taken to be the set of resources identified by the
capacity auction and where the $x$- and $y$-axes give respectively the
powers of the additional stores~$i$ and~$j$.  These additional stores
all have capacities of 100 MWh energy and levels of power which vary
from $10$ to~$50$ MW.  The lowest powered stores contribute only
modest additional capacity, and here~\eqref{eq:3} is seen to be very
accurate.  The highest powered stores contribute substantial
additional capacity; nevertheless here the difference between the left
and right sides of~\eqref{eq:3} is at most~$5\%$.  Similar results
would be obtained if the additional resources~$i$ and~$j$ corresponded
to firm capacity or conventional generation.
The smoothness assumption therefore also appears sufficiently well
satisfied here.

\vspace{-2ex}

\section{Conclusions}
\label{sec:conclusions}

We have given a general theory to show how to value both variable
generation and storage so as to enable them to be integrated fairly
and optimally into electricity capacity markets, whether the objective
within such a market is the optimal satisfaction of some
security-of-supply metric or of some appropriate economic criterion.
We have also considered the consequences for the practical operation
of such markets.
Storage, in particular, provides considerable challenges.  We have
argued that an appropriate risk metric is EEU rather than LOLE, shown
how to calculate the latter for storage, and shown how to value its
economic contribution.  The EFC of storage is sensitive to the amount
of storage already in the system to which it is being added, and this
has considerable consequences for the fair and optimal operation of
markets, as we demonstrate in a realistic practical example based on
the current GB system (in which both conventional generation and
storage---but not currently variable generation---participate).

\vspace{-2ex}

\section*{Acknowledgements}
\label{sec:acknowledgement}

The authors would like to thank the Isaac Newton Institute for
Mathematical Sciences for support during the programme Mathematics of
Energy Systems, when work on this paper was undertaken.  This work was
supported by EPSRC grants numbers EP/I017054/1, EP/K002252/1,
EP/R014604/1, EP/N030028/1 and EP/P001173/1 and by National Grid ESO.
The authors are grateful to many colleagues---David Newbery, Benjamin
Hobbs, Muireann Lynch, Daniel Burke and colleagues at National Grid
ESO---for helpful comments and discussions.  Finally, they are most
grateful to the reviewers for many insightful comments and suggestions
for improvements.

\vspace{-2ex}




\printbibliography


\newpage

\appendix


\section*{Appendix: technical results}

In this appendix we formalise and prove two technical results in the
present paper.




\vspace{-1 ex}

\paragraph{Proof of equation~\eqref{eq:9}.}
\label{sec:proof-equat-eqref}

We prove the result given by equation~\eqref{eq:9}, namely that for
any set of capacity-providing resources~$\R$ which consists entirely
of generation, either conventional or variable, we have that the
derivative $\eeu'(\R)$ of $\eeu(\R)$ with respect to variation of firm
capacity is given by $-\lole(\R)$.

For each time~$t$ let the random variable~$Z_t$ be the supply-demand
balance at time~$t$ corresponding to the use of the set of
resources~$\R$.   Then, from~\eqref{eq:2}, for any addition to $\R$ of
firm capacity equal to $\delta$,
\begin{displaymath}
  \eeu(\R+\delta)
  = \sum_{t=1}^n \int_{-\infty}^0 \Pr(Z_t + \delta < z)\,dz
  = \sum_{t=1}^n \int_{-\infty}^{-\delta} \Pr(Z_t < z)\,dz,
\end{displaymath}
so that, differentiating with respect to~$\delta$ and then setting
$\delta=0$, we have
\begin{displaymath}
  \eeu'(\R) = -\sum_{t=1}^n \Pr(Z_t < 0) = -\lole(\R)
\end{displaymath}
from~\eqref{eq:1} as required.

\paragraph{Proof of Proposition~\ref{prop:one-one}.}

For any possible \emph{risk level}~$k$ of EEU, define~$\R_k$ to be the
set of resources~$\R$ which minimises the cost~$c(\R)$ subject to the
constraint $\eeu(\R)=k$.   Given the risk level~$k$, the subproblem of
determining~$\R_k$ is the problem considered in
Section~\ref{sec:capacity-markets}.   The additional problem in the
minimisation of the overall economic cost~\eqref{eq:13} may therefore
be viewed as that of determining the value of~$k$ such that~$\R_k$
minimises
\begin{equation}
  \label{eq:16}
  \voll\times\eeu(\R_k) + c(\R_k),
\end{equation}
(with $\eeu(\R_k)=k$) i.e.\ that of determining the optimal level of
EEU to be then obtained at minimum cost.
We may consider the effect of infinitesimal variation of the risk
level~$k$ by considering the corresponding required infinitesimal
variation in EFC, where the latter is defined with respect to~$\R_k$.
At that value of~$k$ such that the overall economic cost~\eqref{eq:16}
is minimised, we have stationarity with respect to such variation, so
that at this value of~$k$, analogously to~\eqref{eq:14},
\begin{equation}
  \label{eq:17}
  \voll\times\eeu'(\R_k) + c'(\R_k) = 0,
\end{equation}
where it follows from the definition of EFC in
Section~\ref{sec:equiv-firm-capac} that $\eeu'(\R_k)$ may continue to
be interpreted as defined in that section, i.e.\ as the derivative of
$\eeu(\R_k)$ with respect to firm capacity, and where~$c'(\R_k)$ may
similarly continue to be interpreted as the cost of new entry
($\cone$) at the level of resource defined by~$\R_k$.   Since it is
assumed that all capacity-providing resource consists of some form of
\emph{generation}, the result~\eqref{eq:9} holds and so, analogously
to~\eqref{eq:15}, it follows that at the value of~$k$ such that the
overall cost~\eqref{eq:13} is minimised,
\begin{equation}
  \label{eq:18}
  \lole(\R_k) = \frac{\cone}{\voll}.
\end{equation}

Now, for each~$k$, recall the above definition of $\R_k$.  It follows
from the assumed one-one correspondence between values of $\eeu(\R)$
and those of $\lole(\R)$ that $\R_k$ is also the set of resources~$\R$
which minimises~$c(\R)$ subject to the corresponding constraint on
$\lole(\R)$.  It thus follows from~\eqref{eq:18} that, exactly as in
the case where all resource is provided by firm capacity only, the
determination of the optimal set of capacity-providing resources
minimising the overall economic cost~\eqref{eq:13} is again given by
the minimisation of the cost~$c(\R)$ subject to the
constraint~\eqref{eq:15}.


\end{document}